\documentclass[12pt]{article}
\usepackage[pctex32]{graphics}
\usepackage{amsthm,amsmath,amssymb,amscd,verbatim,epsfig}
\usepackage{amssymb}
\usepackage{graphicx}
\usepackage{epsfig}
\usepackage{graphics,pifont}
\newcommand{\beq}{\begin{equation}}
\newcommand{\eeq}{\end{equation}}

\date{}

\newcommand{\f}{\frac}
\newcommand{\ra}{\rightarrow}

\begin{document}

\title{A\ Simple\ Method\ Which\ Generates\ Infinitely\ Many\ Congruence\ Identities}
\author{Bau-Sen Du \\ [.5cm]
Institute of Mathematics \\
Academia Sinica \\
Taipei 11529, Taiwan \\
dubs@math.sinica.edu.tw \\ [.2cm]
(Fibonacci Quarterly 27 (1989), 116-124) \\}
\maketitle
\begin{abstract}
A simple method called symbolic representation for piecewise linear functions on the real line is introduced and used to compute the numbers of periodic points of all periods for some such functions.  Since, for every positive integer $m$, the number of periodic points of minimal period $m$ must be divisible by $m$, we obtain infinitely many congruence identities.  
\end{abstract}


\section{Introduction}

Let $\phi(m)$ be an integer-valued function defined on the set of all positive integers.  If $m=p_1^{k_1}p_2^{k_2} \cdots p_r^{k_r}$, where the $p_i$'s are distinct prime numbers, $r$ and the $k_i$'s are positive integers, we define $\Phi_1(1, \phi)=\phi(1)$ and $\Phi_1(m, \phi) =$
$$
\phi(m)-\sum_{i=1}^r \phi(\f m{p_i})+\sum_{i_1<i_2} \phi(\f m{p_{i_1}p_{i_2}})
- \sum_{i_1<i_2<i_3} \phi(\f m{p_{i_1}p_{i_2}p_{i_3}}) + \cdots 
+ (-1)^r \phi(\f m{p_1p_2 \cdots p_r}),
$$
\noindent
where the summation $\sum_{i_1<i_2< \cdots < i_j}$ is taken over all integers $i_1, i_2, \cdots, i_j$ with $1 \le i_1 < i_2 <$ $\cdots < i_j \le r$.  

If $m = 2^{k_0}p_1^{k_1}p_2^{k_2} \cdots p_r^{k_r}$, where the $p_i$'s are distinct odd prime numbers, and $k_0 \ge 0, r\ge
1$, and the $k_i$'s $\ge 1$ are integers, we define, similarly, $\Phi_2(m, \phi) =$
$$
\phi(m)-\sum_{i=1}^r \phi(\f m{p_i})+\sum_{i_1<i_2} \phi(\f m{p_{i_1}p_{i_2}})
- \sum_{i_1<i_2<i_3} \phi(\f m{p_{i_1}p_{i_2}p_{i_3}}) + \cdots 
+ (-1)^r \phi(\f m{p_1p_2 \cdots p_r}), 
$$

\noindent
If $m = 2^k$, where $k \ge 0$ is an integer, we define $\Phi_2(m, \phi) = \phi(m) - 1$.  

\noindent
If, for some integer $n \ge 2$, we have $\phi(m) = n^m$ for all positive integer $m$, then we denote $\Phi_i(m, \phi)$ by $\Phi_i(m, n)$, $i = 1, 2$ to emphasize the role of this integer $n$.  

On the other hand, let $S$ be a subset of the real numbers and let $f$ be a function from $S$ into itself.  For every positive integer $n$, we let $f^n$ denote the $n^{\text{th}}$ iterate of $f$: $f^1 = f$ and $f^n = f \circ f^{n-1}$ for $n \ge 2$.  For every $x_0 \in S$, we call the set $\{ f^k(x_0) : k \ge 0 \}$ the orbit of $x_0$ under $f$.  If $x_0$ satisfies $f^m(x_0) = x_0$ for some positive integer $m$, then we call $x_0$ a periodic point of $f$ and call the smallest such positive integer $m$ the minimal period of $x_0$ and of the orbit of $x_0$ (under $f$).  Note that, if $x_0$ is a periodic point of $f$ with minimal period $m$, then, for every integer $1 \le k \le m$, $f^k(x_0)$ is also a periodic point of $f$ with minimal period $m$ and they are all distinct.  So, every periodic orbit of $f$ with minimal period $m$ consists of exactly $m$ distinct points.  Since it is obvious that distinct periodic orbits of $f$ are pairwise disjoint, the number (if finite) of distinct periodic points of $f$ with minimal period $m$ is divisible by $m$ and the quotient equals the number of distinct periodic orbits of $f$ with minimal period $m$.  This observation, together with a standard inclusion-exclusion argument, gives the following well-known result.  
                                                     
\noindent
{\bf Theorem 1.}
{\it Let $S$ be a subset of the real numbers and let $f : S \ra S$ be a mapping with the property that, for every positive integer $m$, the equation $f^m(x) = x$ \, (or $-x$, respectively) has only finitely many distinct solutions.  Let $\phi(m)$ (or $\psi(m)$, respectively) denote the number of these solutions.  Then, for every positive integer $m$, the following hold:
\begin{itemize}

\item[(i)] 
The number of periodic points of $f$ with minimal period $m$ is $\Phi_1(m, \phi)$.  So $\Phi_1(m, \phi) \equiv 0 \, (\text{mod} \,\,\, m)$.

\item[(ii)]
If $0 \in S$ and $f$ is odd, then the number of symmetric periodic points (i.e., periodic points whose orbit are symmetric with respect to the origin) of $f$ with minimal period $2m$ is $\Phi_2(m, \psi)$.  Thus, $\Phi_2(m, \psi) \equiv 0 \, (\text{mod} \,\,\, 2m)$.
\end{itemize}}

Successful applications of the above theorem depend of course on a knowledge of the function $\phi$ or $\psi$.  For example, if we let $S$ denote the set of all real numbers and, for every integer $n \ge 2$ and every odd integer $t = 2k+1 > 1$, let $$f_n(x) = a_n \cdot \Pi_{j=1}^n (x-j)$$ and let $$g_t(x) = b_t \cdot x \cdot \Pi_{j=1}^k (x^2-j^2),$$ where $a_n$ and $b_t$ are fixed sufficiently large positive numbers depending only on $n$ and $t$, respectively.  Then it is easy to see that, for every positive integer $m$, the equation $f_n^m(x) = x \,\,\, (g_t^m(x) = -x$, resp.) has exactly $n^m \,\,\, (t^m$, resp.) distinct solutions in $S$.  Therefore, if $\phi(m,n) = n^m$ and $\psi(m,t) = t^m$, then we have as a consequence of Theorem 1 the following well-known congruence identities which include Fermat's Little Theorem as a special case.

\noindent
{\bf Corollary 2.} {\it (i) Let $m \ge 1$ and $n \ge 2$ be integers.  Then $\Phi_1(m, n) \equiv 0 \,\,(\text{mod} \,\,\, m)$.

\qquad\qquad\,\, (ii) Let $m \ge 1$ be an integer and let $n > 1$ be an odd integer.  Then $\Phi_2(m, n) \equiv 0 \,\,(\text{mod} \,\,\,2m)$.}

In this note, we indicate that the method introduced in {\bf{\cite{du1}}} can also be used to recursively define infinitely many $\phi$ and $\psi$ and thus produce infinitely many families of congruence identities related to Theorem 1.  In Section 2, we will review this method, and to illustrate it we will prove the following result in Section 3.  

\noindent
{\bf Theorem 3.}
{\it For every positive integer $n \ge 3$, let $\phi_n$ be the integer-valued function on the set of all positive integers defined recursively by letting $\phi_n(m) = 2^m - 1$ for all $1 \le m \le n-1$ and $$\phi_n(n+k) = \sum_{j=1}^{n-1} \phi_n(n+k-j), \,\,\, \text{for all} \,\,\, k \ge 0.$$  Then, for every positive integer $m$, $\Phi_1(m, \phi_n) \equiv 0 \,\,\, (\text{mod} \,\,\, m)$.  Furthermore, $$\lim_{m \to \infty} [\log \Phi_1(m, \phi_n)]/m = \lim_{m \to \infty} [\log \phi_n(m)]/m = \log \alpha_n,$$ where $\alpha_n$ is the (unique) positive (and the largest in absolute value) zero of the polynomial $$x^{n-1} - \sum_{k=0}^{n-2} x^k.$$}

Note that in the above theorem these numbers $\phi_n(m)$, $m \ge 1$ are generalized Fibonacci numbers {\bf{\cite{ga,mi}}} and when $n = 3$, these numbers $\phi_3(m)$, $m \ge 1$, are the well-known Lucas numbers: $1, 3, 4, 7, 11, 18, 29, \cdots$.

Just for comparison, we also include the following two results which can be verified numerically.  The rigorous proofs of these two results which are similar to that of Theorem 3 below can be found in ${\bf [1}$, Theorem 2 ${\bf ]}$ and ${\bf [2}$, Theorem 3 ${\bf ]}$, respectively.

\noindent
{\bf Theorem 4.}
{\it For every positive integer $n \ge 2$, let sequences $$<b_{k,1,j,n}>, <b_{k,2,j,n}>, 1\le j \le n,$$ be defined recursively as follows:
$$\begin{cases}
b_{1,1,j,n} =  0, \,\,\, 1 \le j \le n, \\
b_{2,1,j,n} =  1, \,\,\, 1 \le j \le n, \\
b_{1,2,j,n} =  b_{2,2,j,n} = 0, \,\,\, 1 \le j \le n-1, \\
b_{1,2,n,n} =  b_{2,2,n,n} = 1. \\
\end{cases}
$$  For $i = 1$ or 2, and $k \ge 1$, 
$$\begin{cases}
b_{k+2,i,j,n} = b_{k,i,1,n} + b_{k,i,j+1,n}, \,\,\, 1 \le j \le n-1, \\
b_{k+2,i,n,n} = b_{k,i,1,n} + b_{k+1,i,n,n}. \\
\end{cases}
$$  Let $b_{k,1,j,n} = 0$ for all $-2n+3 \le k \le 0$ and $1 \le j \le n$, and for all positive integers $m$, let $$\phi_n(m) = b_{m,2,n,n} + 2 \cdot \sum_{j=1}^n b_{m+2-2j,1,j,n}.$$ Then, for every positive integer $m$, $\Phi_1(m, \phi_n) \equiv 0 \,\,\, (\text{mod} \,\,\, m)$.  Furthermore, $$\lim_{m \to \infty} [\log \Phi_1(m, \phi_n)]/m = \lim_{m \to \infty} [\log \phi_n(m)]/m = \log \beta_n,$$ where $\beta_n$ is the (unique) positive (and the largest in absolute value) zero of the polynomial $x^{2n+1} - 2x^{2n-1} - 1$.}

\noindent
{\bf Remark 1.}
For all positive integers $m$ and $n$, let $$A_{m,n} = \Phi_1(2m-1, \phi_n)/(2m-1),$$ where $\phi_n$ is defined as in Theorem 3 for $n = 1$ and as in Theorem 4 for $n \ge 2$.  Table 1 lists the first 31 values of $A_{m,n}$, for $1 \le n \le 6$.  It seems that $A_{m,n} = 2^{m-n-1}$ for $n+1 \le m \le 3n+2$ and $A_{m,n} > 2^{m-n-1}$ for $m > 3n + 2$.  If, for all positive integers $m$ and $n$, we define sequences $<B_{m,n,k}>$ by letting $$B_{m,n,1} = A_{m+3n+2,n} - 2A_{m+3n+1,n} \qquad\quad\,\,$$ and $$B_{m,n,k} = B_{m+2n+1,n,k-1} - B_{m+2n+1,n+1,k-1}$$ for $k > 1$, then more extensive numerical computations seem to show that, for all positive integers $k$, we have

\noindent
(i) \,\, $B_{1,n,k} = 2$ \,\,\,\,\,\, for all \,\, $n \ge 1$, \\
(ii) \, $B_{2,n,k} = 4k$ \,\, for all \,\, $n \ge 1$, \\
(iii) $B_{3,n,k}$ is a constant depending only on $k$, and \\
(iv) for all \,\, $1 \le m \le 2n+1$, \,\, $B_{m,n,k} = B_{m,j,k}$ \,\, for all \,\, $j \ge n \ge 1$. \\

\noindent
{\bf Theorem 5.}
{\it Fix any integer $n \ge 2$.  For all integers $i$, $j$, and $k$ with $i = 1, 2$, $1 \le |j| \le n$, and $k \ge 1$, we define $c_{k,i,j,n}$ recursively as follows: $$c_{1,1,n,n} = 1 \,\,\, \text{and} \,\,\, c_{1,1,j,n} = 0 \,\,\, \text{for} \,\,\, j \ne n,$$ $$c_{1,2,1,n} = 1 \,\,\, \text{and} \,\,\, c_{1,2,j,n} = 0 \,\,\, \text{for} \,\,\, j \ne 1,$$  For $i = 1, 2$, and $k \ge 1$, 
$$
\begin{cases}
c_{k+1,i,1,n} &= c_{k,i,1,n} + c_{k,i,-n,n} + c_{k,i,n,n}, \\
c_{k+1,i,j,n} &= c_{k,i,j-1,n} + c_{k,i,n,n} \,\,\, \text{for all} \,\,\, 2\le j \le n, \\
c_{k+1,i,-1,n} &= c_{k,i,-1,n} + c_{k,i,-n,n} + c_{k,i,n,n}, \\
c_{k+1,i,-j,n} &= c_{k,i,-j+1,n} + c_{k,i,-n,n} \,\,\, \text{for all} \,\,\, 2 \le j \le n. \\
\end{cases}
$$  Let $c_{k,1,j,n} = 0$ for all integers $k$, $j$ with $4-n \le k \le 0$ and $1 \le |j| \le n$, and, for all positive integers $m$, let $$ \phi_n(m) = 2 \sum_{k=1}^{n-1} c_{m+2-k,1,n+1-k,n} + 2c_{m+1,2,1,n} - 1\,\,\,\,$$ and $$\psi_n(m) = 2 \sum_{k=1}^{n-1} c_{m+2-k,1,k-n-1,n} + 2c_{m+1,2,-1,n} + 1.$$  Then, for every positive integer $m$, $$\Phi_1(m, \phi_n) \equiv 0 \,\,\, (\text{mod} \,\,\, m) \quad \text{and} \quad \Phi_2(m, \psi_n) \equiv 0 \,\,\, (\text{mod} \,\,\, 2m).$$  Furthermore, $$\lim_{m \to \infty} [\log \Phi_1(m, \phi_n)]/m = \lim_{m \to \infty} [\log \phi_n(m)]/m = \lim_{m \to \infty} [\log \psi_n(m)]/m$$ $$\qquad\qquad\qquad= \lim_{m \to \infty} [\log \Phi_2(m, \psi_n)]/m = \log \gamma_n.$$ where $\gamma_n$ is the (unique) positive (and the largest in absolute value) zero of the polynomial $x^n - 2x^{n-1} - 1$.} 

\noindent
{\bf Remark 2.}
For all positive integers $m \ge 1$ and $n \ge 2$, let $$D_{m,n} = \Phi_2(m, \psi_n)/(2m),$$ where the $\psi_n$'s are defined as in the above theorem.  Table 2 lists the first 25 values of $D_{m,n}$ for $2 \le n \le 6$.  It seems that $D_{m,n} = 2^{m-n}$ for $n \le m \le 3n$, and $D_{m,n} > 2^{m-n}$ for $m > 3n$.  If, for all integers $m \ge 1$ and $n \ge 2$, we define the sequences $<E_{m,n,k}>$ by letting $$E_{m,n,1} = D_{m+3n,n} - 2D_{m+3n-1,n}\qquad\,\,$$ and $$E_{m,n,k} = E_{m+2n,n,k-1} - E_{m+2n,n+1,k-1}$$ for $k > 1$, then more extensive computations seem to show that, for all positive integers $k$, we have

\noindent
(i) \,\, $E_{1,n,k} = 2$ \,\,\,\,\,\, for all \,\, $n \ge 2$, \\
(ii) \, $E_{2,n,k} = 4k$ \,\, for all \,\, $n \ge 2$, \\
(iii) $E_{3,n,k}$ and $E_{4,n,k}$ are constants depending only on $k$, and \\
(iv) for all \,\, $1 \le m \le 2n$, \,\, $E_{m,n,k} = E_{m,j,k}$ \,\, for all \,\, $j \ge n \ge 2$. \\

\noindent
See Tables 1 and 2 below.
\pagebreak

$$\text{Table 1}$$
\bigskip
\noindent
$m$ \qquad\qquad $A_{m,1}$ \qquad\qquad $A_{m,2}$ \qquad\qquad $A_{m,3}$ \qquad\qquad $A_{m,4}$ \qquad\qquad $A_{m,5}$ \qquad\qquad $A_{m,6}$

\noindent
1 \quad \qquad\qquad 1 \qquad\qquad\qquad 1 \quad\qquad\qquad 1 \quad\qquad\qquad 1 \quad\quad\qquad\qquad 1 \quad\qquad\qquad 1 \\
2 \quad\qquad\qquad 1 \qquad\qquad\qquad 1 \quad\qquad\qquad 1 \quad\qquad\qquad 1 \quad\quad\qquad\qquad 1 \quad\qquad\qquad 1 \\  
3 \quad\qquad\qquad 1 \qquad\qquad\qquad 0 \quad\qquad\qquad 0 \quad\qquad\qquad 0 \quad\quad\qquad\qquad 0 \quad\qquad\qquad 0 \\
4 \quad\qquad\qquad 1 \qquad\qquad\qquad 1 \quad\qquad\qquad 1 \quad\qquad\qquad 1 \quad\quad\qquad\qquad 1 \quad\qquad\qquad 1 \\
5 \quad\qquad\qquad 2 \qquad\qquad\qquad 1 \quad\qquad\qquad 0 \quad\qquad\qquad 0 \quad\quad\qquad\qquad 0 \quad\qquad\qquad 0 \\
6 \quad\qquad\qquad 2 \qquad\qquad\qquad 2 \quad\qquad\qquad 2 \quad\qquad\qquad 2 \quad\quad\qquad\qquad 2 \quad\qquad\qquad 2 \\
7 \quad\qquad\qquad 4 \qquad\qquad\qquad 2 \quad\qquad\qquad 1 \quad\qquad\qquad 0 \quad\quad\qquad\qquad 0 \quad\qquad\qquad 0 \\
8 \quad\qquad\qquad 5 \qquad\qquad\qquad 3 \quad\qquad\qquad 3 \quad\qquad\qquad 3 \quad\quad\qquad\qquad 3 \quad\qquad\qquad 3 \\
9 \quad\qquad\qquad 8 \qquad\qquad\qquad 4 \quad\qquad\qquad 2 \quad\qquad\qquad 1 \quad\quad\qquad\qquad 0 \quad\qquad\qquad 0 \\
10 \qquad\quad\quad\,\, 11 \qquad\qquad\qquad\, 6 \quad\quad\qquad 6 \quad\qquad\qquad 6 \quad\quad\qquad\qquad 6 \quad\qquad\qquad 6 \\
11 \qquad\quad\quad\,\, 18 \qquad\qquad\qquad\, 8 \quad\quad\qquad 4 \quad\qquad\qquad 2 \quad\quad\qquad\qquad 1 \quad\qquad\qquad 0 \\
12 \qquad\quad\quad\,\, 25 \qquad\qquad\qquad 11 \quad\quad\qquad 9 \quad\qquad\qquad 9 \quad\quad\qquad\qquad 9 \,\,\qquad\qquad\,\, 9 \\
13 \qquad\quad\quad\,\, 40 \qquad\qquad\qquad 16 \quad\quad\qquad 8 \quad\qquad\qquad 4 \quad\quad\qquad\qquad 2 \,\,\qquad\qquad\,\, 1 \\
14 \qquad\quad\quad\,\, 58 \qquad\qquad\qquad 23 \quad\quad\qquad 18 \qquad\qquad 18 \quad\quad\qquad\qquad 18 \qquad\qquad 18 \\
15 \qquad\quad\quad\,\, 90 \qquad\qquad\qquad 32 \quad\quad\qquad 16 \quad\qquad\qquad 8 \quad\quad\qquad\qquad 4 \qquad\qquad\,\, 2 \\
16 \qquad\quad\quad 135 \qquad\qquad\qquad 46 \quad\quad\qquad 32 \quad\qquad\qquad 30 \quad\qquad\qquad 30 \qquad\qquad 30 \\
17 \qquad\quad\quad 210 \qquad\qquad\qquad 66 \quad\quad\qquad 32 \quad\qquad\qquad 16 \quad\qquad\qquad 8 \qquad\qquad\,\,\,\, 4 \\
18 \qquad\quad\quad 316 \qquad\qquad\qquad 94 \quad\quad\qquad 61 \quad\qquad\qquad 56 \quad\qquad\qquad 56 \qquad\qquad 56 \\
19 \qquad\quad\quad 492 \quad\qquad\qquad\,\, 136 \quad\quad\qquad 64 \quad\qquad\qquad 32 \quad\qquad\qquad 16 \qquad\qquad\,\, 8 \\
20 \qquad\quad\quad 750 \quad\qquad\qquad\,\, 195 \quad\qquad\,\,\, 115 \quad\qquad\qquad 101 \quad\qquad\qquad 99 \qquad\qquad 99 \\
21 \,\,\,\quad\quad\quad 1164 \quad\qquad\qquad\,\, 282 \quad\quad\qquad 128 \quad\qquad\qquad 64 \quad\qquad\qquad 32 \qquad\qquad 16 \\
22 \,\,\quad\quad\quad 1791 \quad\qquad\qquad\,\, 408 \quad\quad\qquad 224 \qquad\qquad\,\, 191 \quad\qquad\qquad 186 \quad\qquad\,\, 186 \\
23 \,\,\quad\quad\quad 2786 \quad\qquad\qquad\,\, 592 \quad\quad\qquad 258 \qquad\qquad\,\, 128 \quad\qquad\qquad\,\, 64 \qquad\qquad 32 \\
24 \,\,\quad\quad\quad 4305 \quad\qquad\qquad\,\, 856 \quad\quad\qquad 431 \qquad\qquad\,\, 351 \quad\qquad\qquad 337 \quad\qquad\,\, 335 \\
25 \,\,\quad\quad\quad 6710 \quad\qquad\qquad 1248 \qquad\qquad\,\, 520 \qquad\qquad 256 \quad\qquad\qquad 128 \qquad\qquad 64 \\
26 \quad\quad\quad 10420 \,\,\,\qquad\qquad\,\, 1814 \quad\qquad\,\,\,\, 850 \qquad\qquad\,\, 668 \qquad\qquad\,\, 635 \qquad\qquad 630 \\
27 \quad\quad\quad 16264 \quad\qquad\qquad 2646 \qquad\qquad 1050 \qquad\qquad\,\, 512 \qquad\qquad\,\, 256 \qquad\qquad 128 \\
28 \quad\quad\quad 25350 \quad\qquad\qquad 3858 \qquad\qquad 1673 \qquad\qquad 1257 \qquad\qquad 1177 \qquad\qquad 1163 \\
29 \quad\quad\quad 39650 \quad\qquad\qquad 5644 \qquad\qquad 2128 \qquad\qquad 1026 \qquad\qquad\,\, 512 \qquad\qquad 256 \\
30 \quad\quad\quad 61967 \quad\qquad\qquad 8246 \qquad\qquad 3328 \qquad\qquad 2402 \qquad\qquad 2220 \qquad\qquad 2187 \\
31 \quad\quad\quad 97108 \qquad\qquad\,\, 12088 \qquad\qquad 4320 \qquad\qquad 2056 \qquad\qquad 1024 \qquad\qquad 512 \\
\pagebreak

$$\text{Table 2}$$
\bigskip
\noindent
$m$ \quad\quad\qquad\quad $D_{m,2}$ \qquad\qquad\qquad\,\, $D_{m,3}$ \qquad\qquad\qquad $D_{m,4}$ \quad\qquad\qquad $D_{m,5}$ \quad\qquad\qquad $D_{m,6}$

\noindent
1 \quad\qquad\qquad\qquad\,\,\,\, 0 \quad\qquad\qquad\qquad\,\,\,\,\,\, 0 \qquad\qquad\qquad\qquad 0 \qquad\quad\quad\quad\,\, 0 \quad\quad\qquad\qquad 0  \\
2 \quad\qquad\qquad\qquad\,\,\,\, 1 \quad\qquad\qquad\qquad\,\,\,\,\,\, 0 \qquad\qquad\qquad\qquad 0 \qquad\quad\quad\quad\,\, 0 \quad\quad\qquad\qquad 0 \\  
3 \quad\qquad\qquad\qquad\,\,\,\, 2 \quad\qquad\qquad\qquad\,\,\,\,\,\, 1 \qquad\qquad\qquad\qquad 0 \qquad\quad\quad\quad\,\, 0 \quad\quad\qquad\qquad 0 \\
4 \quad\qquad\qquad\qquad\,\,\,\, 4 \quad\qquad\qquad\qquad\,\,\,\,\,\, 2 \qquad\qquad\qquad\qquad 1 \qquad\quad\quad\quad\,\, 0 \quad\quad\qquad\qquad 0  \\
5 \quad\qquad\qquad\qquad\,\,\,\, 8 \quad\qquad\qquad\qquad\,\,\,\,\,\, 4 \qquad\qquad\qquad\qquad 2 \qquad\quad\quad\quad\,\, 1 \quad\quad\qquad\qquad 0  \\
6 \quad\qquad\qquad\qquad\, 16 \,\,\,\qquad\qquad\qquad\quad\,\,\, 8 \qquad\qquad\qquad\qquad 4 \qquad\quad\qquad\quad\,\, 2 \quad\quad\quad\qquad 1  \\
7 \quad\qquad\qquad\qquad\, 34 \quad\quad\qquad\qquad\quad\,\,\, 16 \qquad\qquad\qquad\qquad 8 \qquad\quad\qquad\quad\,\, 4 \quad\quad\quad\qquad 2  \\
8 \quad\qquad\qquad\qquad\,\, 72 \quad\qquad\qquad\qquad\,\, 32 \quad\qquad\qquad\qquad\,\, 16 \qquad\quad\qquad\quad\,\, 8 \quad\quad\quad\qquad 4  \\
9 \quad\qquad\qquad\qquad 154 \quad\qquad\qquad\qquad\,\, 64 \quad\qquad\qquad\qquad\,\, 32 \qquad\quad\qquad\quad 16 \quad\quad\quad\qquad 8  \\
10 \qquad\quad\quad\qquad\,\, 336 \quad\qquad\qquad\qquad 130 \quad\qquad\qquad\qquad\,\, 64 \qquad\quad\qquad\quad 32 \quad\quad\quad\qquad 16  \\
11 \qquad\quad\quad\qquad\,\, 738 \qquad\qquad\qquad\quad 264 \quad\qquad\qquad\qquad 128 \quad\qquad\qquad\,\,\,\,\, 64 \quad\quad\quad\qquad 32  \\
12 \qquad\quad\quad\qquad 1632 \qquad\qquad\qquad\quad 538 \quad\qquad\qquad\qquad 256 \quad\qquad\qquad\,\,\, 128 \quad\quad\quad\qquad 64  \\
13 \qquad\quad\quad\qquad 3640 \qquad\qquad\qquad\,\, 1104 \quad\qquad\qquad\qquad 514 \quad\qquad\qquad\,\,\, 256 \quad\quad\quad\quad 128  \\
14 \qquad\quad\quad\qquad 8160 \qquad\qquad\qquad\,\, 2272 \quad\quad\qquad\qquad\, 1032 \qquad\qquad\qquad 512 \quad\quad\quad\quad 256 \\
15 \qquad\quad\quad\quad\,\, 18384 \qquad\qquad\qquad\,\, 4692 \quad\quad\qquad\qquad\, 2074 \quad\qquad\qquad\,\, 1024 \quad\quad\quad\quad 512  \\
16 \qquad\quad\quad\quad\, 41616 \qquad\qquad\qquad\,\, 9730 \qquad\qquad\qquad\,\, 4176 \quad\qquad\qquad\,\, 2050 \quad\quad\qquad 1024 \\
17 \qquad\quad\quad\,\,\quad 94560 \qquad\qquad\qquad 20236 \qquad\qquad\qquad\,\, 8416 \quad\qquad\qquad\,\, 4104 \quad\quad\qquad 2048 \\
18 \qquad\quad\quad\quad 215600 \qquad\qquad\qquad 42208 \qquad\qquad\qquad 16980 \quad\qquad\qquad\,\, 8218 \quad\quad\qquad 4096 \\
19 \qquad\quad\quad\,\,\,\,\, 493122 \qquad\qquad\qquad 88288 \qquad\qquad\qquad 34304 \quad\qquad\qquad 16464 \quad\quad\qquad 8194 \\
20 \qquad\quad\quad\,\, 1130976 \qquad\qquad\quad\,\, 185126 \qquad\qquad\qquad 69376 \quad\qquad\qquad 32992 \quad\quad\qquad 16392  \\
21 \qquad\quad\quad\,\, 2600388 \quad\qquad\qquad\,\, 389072 \quad\qquad\qquad\,\, 140458 \quad\qquad\qquad 66132 \quad\quad\qquad 32794  \\
22 \qquad\quad\quad\,\, 5992560 \quad\qquad\qquad\,\, 819458 \quad\qquad\qquad\,\, 284684 \qquad\qquad\,\, 132608 \quad\quad\qquad 65616 \\
23 \qquad\quad\quad 13838306 \quad\qquad\qquad 1729296 \quad\qquad\qquad\,\, 577592 \qquad\qquad\,\, 265984 \quad\quad\qquad 131296  \\
24 \qquad\quad\quad 32016576 \quad\qquad\qquad 3655936 \quad\qquad\qquad 1173040 \qquad\qquad\,\, 533672 \quad\quad\qquad 262740 \\
25 \qquad\quad\quad 74203112 \quad\qquad\qquad 7742124 \quad\qquad\qquad 2384678 \qquad\qquad 1071104 \quad\quad\qquad 525824 \\

\section{Symbolic representation for continuous piecewise linear functions}

In this section, we review the method introduced in {\bf{\cite{du1}}}.  Throughout this section, let $g$ be a continuous piecewise linear function from the interval $[c, d]$ into itself.  We call the set $\{(x_i, y_i) : i = 1,2, \cdots, k \}$ a set of nodes for (the graph of) $y = g(x)$ if the following three conditions hold:
\begin{itemize}

\item[(1)]
$k \ge 2$,

\item[(2)]
$x_1 = c$, $x_k = d$, $x_1 < x_2 < \cdots < x_k$, and 

\item[(3)]
$g$ is linear on $[x_i, x_{i+1}]$ for all $1 \le i \le k-1$ and $y_i = g(x_i)$ for all $1 \le i \le k$.
\end{itemize}

\noindent
For any such set, we will use its $y$-coordinates $y_1, y_2, \cdots, y_k$ to represent its graph and call $y_1y_2 \cdots y_k$ (in that order) a (symbolic) representation for (the graph of) $y = g(x)$.  For $1 \le i < j \le k$, we call $y_iy_{i+1} \cdots y_j$ the representation for $y = g(x)$ on $[x_i, x_j]$ obtained by restricting $y_1y_2 \cdots y_k$ to $[x_i, x_j]$.  For convenience, we will also call every $y_i$ in $y_1y_2 \cdots y_k$ a node.  If $y_i = y_{i+1}$ for some $i$ (i.e., $g$ is constant on $[x_i, x_{i+1}]$), we will simply write $$y_1 \cdots y_iy_{i+2} \cdots y_k$$ instead of 
$$y_1 \cdots y_iy_{i+1}y_{i+2} \cdots y_k.$$  That is, we will delete $y_{i+1}$ from the (symbolic) representation $y_1y_2 \cdots y_k$.  Therefore, every two consecutive nodes in a (symbolic) representation are distinct.  Note that a continuous piecewise linear function obviously has more than one (symbolic) representation.  However, as we will soon see that there is no need to worry about that.  

Now assume that $\{ (x_i, y_i) : i = 1, 2, \cdots, k \}$ is a set of nodes for $y = g(x)$ and $a_1a_2 \cdots a_r$ is a representation for $y = g(x)$ with $$\{ a_1, a_2, \cdots, a_r \} \subset \{ y_1, y_2, \cdots, y_k \}$$ and $a_i \ne a_{i+1}$ for all $1 \le i \le r-1$.  If $$\{ y_1, y_2, \cdots, y_k \} \subset \{ x_1, x_2, \cdots, x_k \},$$ then there is an easy way to obtain a representation for $y = g^2(x)$ from the one $a_1a_2 \cdots a_r$ for $y = g(x)$.  The procedure is as follows.  First, for any two distinct real numbers $u$ and $v$, let $[u : v]$ denote the closed interval with endpoints $u$ and $v$.  Then let $b_{i,1}b_{i,2} \cdots b_{i,t_i}$ be the representation for $y = g(x)$ on $[a_i:a_{i+1}]$ which is obtained by restricting $a_1a_2 \cdots a_r$ to $[a_i:a_{i+1}]$.  We use the following notation to indicate this fact: $$a_ia_{i+1} \ra b_{i,1}b_{i,2} \cdots b_{i,t_i} \,\,\, (\text{under} \,\,\, g) \,\,\, \text{if} \,\,\, a_i < a_{i+1},$$ or $$a_ia_{i+1} \ra b_{i,t_i} \cdots b_{i,2}b_{i,1} \,\,\, (\text{under} \,\,\, g) \,\,\, \text{if} \,\,\, a_i > a_{i+1}.$$  The above representation on $[a_i:a_{i+1}]$ exists since $$\{ a_1, a_2, \cdots, a_r \} \subset \{ x_1, x_2, \cdots, x_k \}.$$  Finally, if $a_i < a_{i+1}$, let $z_{i,j} = b_{i,j}$ for all $1 \le j \le t_i$.  If $a_i > a_{i+1}$, let $z_{i,j} = b_{i,t_i+1-j}$ for all $1 \le j \le t_i$.  Let $$Z = z_{1,1} \cdots z_{1,t_1}z_{2,2} \cdots z_{2,t_2} \cdots z_{r,2} \cdots z_{r,t_r}.$$  (Note that $z_{i,t_i} = z_{i+1,1}$ \, for all \, $1 \le i \le r-1$).  Then it is easy to see that $Z$ is a representation for $y = g^2(x)$.  It is also obvious that the above procedure can be applied to the representation $Z$ for $y = g^2(x)$ to obtain one for $y = g^3(x)$, and so on.  

\section{Proof of Theorem 3}

In this section we fix an integer $n \ge 3$ and let $f_n(x)$ be the continuous function from the interval $[1, n]$ onto itself defined by $$f_n(x) = x + 1 \,\,\, \text{for} \,\,\, 1 \le x \le n-1 \qquad\qquad\qquad\qquad\,\,\,\,$$ and $$f_n(x) = -(n-1)x + n^2 - n + 1 \,\,\, \text{for} \,\,\, n-1 \le x \le n.$$  Using the notations introduced in Section 2, we have the following result.  

\noindent
{\bf Lemma 6.}
{\it Under $f_n$, we have
$$
\begin{cases}
k(k+1) &\ra (k+1)(k+2), \,\,\, 2 \le k \le n-2, \,\,\, \text{if} \,\,\, n > 3, \\
(k+1)k &\ra (k+2)(k+1), \,\,\, 2 \le k \le n-2, \,\,\, \text{if} \,\,\, n > 3. \\
(n-1)n &\ra n1, \quad n(n-1) \ra 1n, \\
n1 &\ra 1n(n-1) \cdots 432, \quad 1n \ra 234 \cdots (n-1)n1. \\
\end{cases}
$$}

In the following when we say the representation for $y = f_n^k(x)$, we mean the representation obtained, following the procedure as described in Section 2, by applying Lemma 6 to the representation $234 \cdots (n-1)n1$ for $y = f_n(x)$ successively until we get to the one for $y = f_n^k(x)$.

For every positive integer $k$ and all integers $i$, $j$, with $1 \le i, j \le n-1$, let $a_{k,i,j,n}$ denote the number of $uv$'s and $vu$'s in the representation for $y = f_n^k(x)$ whose corresponding $x$-coordinates are in the interval $[i, i+1]$, where $uv = 1n$ if $j = 1$, and $uv = j(j+1)$ if $2 \le j \le n-1$.  It is obvious that $$a_{1,i,i+1,n} = 1 \,\,\, \text{for all} \,\,\, 1 \le i \le n-2, \,\quad\quad$$ $$a_{1,n-1,1,n} = 1, \,\, \text{and} \,\,\, a_{1,i,j,n} = 0 \,\,\, \text{elsewhere}.$$  From the above lemma, we find that these sequences $<a_{k,i,j,n}>$ can be computed recursively.  

\noindent
{\bf Lemma 7.}
{\it For every positive integer $k$ and all integers $i$ with $1 \le i \le n-1$, we have
$$
\begin{cases}
a_{k+1,i,1,n} & = a_{k,i,1,n} + a_{k,i,n-1,n}, \\
a_{k+1,i,2,n} & = a_{k,i,1,n}, \\
a_{k+1,i,j,n} & = a_{k,i,1,n} + a_{k,i,j-1,n}, \,\,\, 3 \le j \le n-1 \,\,\, \text{if} \,\,\, n > 3. \\
\end{cases}
$$}

It then follows from the above lemma that the sequences $<a_{k,i,j,n}>$ can all be computed from the sequences $<a_{k,n-1,j,n}>$.  

\noindent
{\bf Lemma 8.}
{\it For every positive integer $k$ and all integers $j$ with $1 \le j \le n-1$, we have $$a_{k,n-1,j,n} = a_{k+i,n-1-i,j,n}, \,\,\, 1 \le i \le n-2.$$}

For every positive integer $k$, let $$c_{k,n} = \sum_{i=1}^{n-1} a_{k,i,1,n} + \sum_{i=2}^{n-1} a_{k,i,i,n}.$$  Then it is easy to see that $c_{k,n}$ is exactly the number of distinct solutions of the equation $f_n^k(x) = x$ in the interval $[1, n]$.  From the above lemma, we also have, for all $k \ge 1$, the identities: $$c_{k,n} = \sum_{i=0}^{n-2} a_{k-i,n-1,1,n} + \sum_{i=0}^{n-3} a_{k-i,n-1,n-1-i,n}$$ provided that $a_{m,n-1,j,n} = 0$ for all $m \le 0$ and $j > 0$.  Since, for every positive integer $k$, $$a_{k,n-1,1,n} = a_{k-1,n-1,1,n} + a_{k-1,n-1,n-1,n} = a_{k-1,n-1,1,n} + a_{k-2,n-1,1,n} + a_{k-2,n-1,n-2,n}$$ $$= a_{k-1,n-1,1,n} + a_{k-2,n-1,1,n} + a_{k-3,n-1,1,n} + a_{k-3,n-1,n-3,n} = \cdots = \sum_{i=1}^{n-1} a_{k-i,n-1,1,n}$$ and $$c_{k,n} = \sum_{i=0}^{n-2} a_{k-i,n-1,1,n} + \sum_{i=0}^{n-3} a_{k-i,n-1,n-1-i,n}$$ $$= a_{k,n-1,1,n} + a_{k-1,n-1,1,n} + \sum_{i=2}^{n-2} a_{k-i,n-1,1,n} + a_{k-1,n-1,1,n} + a_{k-1,n-1,n-2,n} + \sum_{i=1}^{n-3} a_{k-i,n-1,n-1-i,n}$$ $$= a_{k,n-1,1,n} + 2a_{k-1,n-1,1,n} + \sum_{i=2}^{n-2} a_{k-i,n-1,1,n} + 2a_{k-1,n-1,n-2,n} + \sum_{i=2}^{n-3} a_{k-i,n-1,n-1-i,n} = \cdots$$ $$= \sum_{i=0}^{n-2} (i+1)a_{k-i,n-1,1,n}$$ provided that $a_{m,n-1,1,n} = 0$ if $m \le 0$, we obtain that $c_{k,n} = 2^k - 1$ for all $1 \le k \le n-1$ and $$c_{k,n} = \sum_{i=1}^{n-1} c_{k-i,n} \,\,\, \text{for all integers} \,\,\, k \ge n.$$  If, for every positive integer $m$, we let $\phi_n(m) = c_{m,n}$, then, by Theorem 1, we have $\Phi_1(m, \phi_n) \equiv 0$ (mod $m$).  The proof of the other statement of Theorem 3 is easy and omitted (see {\bf{\cite{ga}}} and {\bf{\cite{mi}}}).  This completes the proof of Theorem 3.

\end{document}